\documentclass[reqno,11pt,a4paper]{amsart}
\usepackage{latexsym,amssymb,amsfonts}
\usepackage{graphicx,color,graphicx}
\usepackage[notref, notcite]{showkeys}
\usepackage{mathrsfs}
\usepackage{lscape}
\newtheorem{lemma}{\bf Lemma}[section]

\newtheorem{remark}{Remark}

\newtheorem{theorem}{\bf Theorem}[section]

\theoremstyle{plain} \numberwithin{equation}{section}
\theoremstyle{plain} \theoremstyle{definition}

\theoremstyle{assumption}
\newtheorem{assumption}{Assumption}[section]
\DeclareMathOperator*{\cov}{cov}
\DeclareMathOperator*{\diag}{diag}
\DeclareMathOperator*{\trace}{trace}
\DeclareMathOperator*{\sign}{sign}

\input cyracc.def

\begin{document}

\title[Moderate deviation principle]
{Moderate deviation principle for ergodic Markov
chain. Lipschitz summands}

\author{B. Delyon}
\address{Universit\'e de Rennes 1, IRISA, Campus de Beaulieu, 35042 Rennes Cedex,
France.} \email{bernard.delyon@univ-rennes1.fr}
\author{A. Juditsky}
\address{University Joseph Fourier of Grenoble, France}
\email{juditsky@inrialpes.fr}

\author{R. Liptser}
\address{Department of Electrical Engineering-Systems,
Tel Aviv University, 69978 Tel Aviv Israel}
\email{liptser@eng.tau.ac.il}

\begin{abstract}
For $\frac{1}{2}<\alpha<1$, we propose the MDP analysis for family
$$
S^\alpha_n=\frac{1}{n^\alpha}\sum_{i=1}^nH(X_{i-1}), \ n\ge 1,
$$
where $(X_n)_{n\ge 0}$ be a homogeneous ergodic Markov chain,
$X_n\in \mathbb{R}^d$,  when the spectrum of operator $P_x$ is
continuous. The vector-valued function $H$ is not assumed to be
bounded but the Lipschitz continuity of $H$ is required. The main
helpful tools in our approach are Poisson's equation and
Stochastic Exponential; the first enables to replace the original
family by $\frac{1}{n^\alpha}M_n$ with a martingale $M_n$ while
the second to avoid the direct Laplace transform analysis.
\end{abstract}
\maketitle
\section{Introduction and discussion}

Let $(X_n)_{n\ge 0}$ be a homogeneous ergodic Markov chain,
$X_n\in \mathbb{R}^d$   with the transition probability kernel for
$n$ steps: $P^{(n)}_x=P^{(n)}(x,dy)$ (for brevity
$P^{(1)}_x:=P_x$) and the unique invariant measure $\mu$.

Let $H$ be a measurable function
$\mathbb{R}^d\stackrel{H}{\to}\mathbb{R}^p$ with
$\int_{\mathbb{R}^d}|H(z)|\mu (dz)<\infty$ and
\begin{equation}\label{nul}
\int_{\mathbb{R}^d}H(z)\mu (dz)=0.
\end{equation}
Set
$$
S^\alpha_n=\frac{1}{n^\alpha}\sum_{i=1}^nH(X_{i-1}), \ n\ge 1; \
(0.5<\alpha<1).
$$

In this paper, we examine the moderate deviation principle (in
short: MDP) for the family $(S^\alpha_n)_{n\ge 1}$ when the
spectrum of operator $P_x$ is continuous.

It is well known that for bounded $H$ satisfying \eqref{nul} ((H)
- condition), the most MDP compatible Markov chains are
characterized by {\it eigenvalues gap condition {\rm (EG)}} (see
Wu, \cite{Wu}, \cite{Wu1}, Gong and Wu, \cite{GongWu5}, and
citations therein):
\begin{quote}
{\it the unit is an isolated, simple and the only eigenvalue with
modulus 1 of the transition probability kernel $P_x$}.
\end{quote}
In the framework of (H)-(EG)  conditions, the MDP is valid with
the rate of speed $n^{-(2\alpha-1)}$ and the rate function
$I(y),y\in\mathbb{R}^d$
\begin{equation}\label{I(x)2}
I(y)=
\begin{cases}
\frac{1}{2}\|y\|^2_{B^\oplus}, & B^\oplus By=y
\\
\infty, & \text{otherwise},
\end{cases}
\end{equation}
where  $B^\oplus$ is the pseudoinverse matrix (in Moore-Penrose
sense, see e.g.\cite {Albert}) for the matrix
\begin{multline}\label{B}
B=\int_{\mathbb{R}^d}H(x)H^*(x)\mu(dx)
\\
+\sum_{n\ge 1}\int_{\mathbb{R}^d}
\left[H(x)(P^{(n)}_xH)^*+(P^{(n)}_xH)H^*(x)\right]\mu(dx)
\end{multline}
(henceforth, $^*$, $|\cdot|$,  and $\|\cdot\|_Q$ are the
transposition symbol, $\mathbb{L}^1$ norm and $\mathbb{L}^2$ norm
with the kernel $Q$ ($\|x\|_Q=\sqrt{\langle x,Qx\rangle}$)
respectively).

Thanks to the quadratic form rate function, the MDP is an
attractive tool for an asymptotic analysis in many areas, say,
with thesis (see, example in Section \ref{sec-7})
\begin{center}
``MDP instead of CLT''.
\end{center}

In this paper, we intend to apply the MDP analysis to Markov chain
defined by the recurrent equation
$$
X_n=f(X_{n-1},\xi_n), \ n\ge 1
$$
generated by i.i.d. sequence $(\xi_n)_{n\ge 1}$ of random vectors,
where $f$ is some vector-valued measurable function. Obviously,
the function $f$ and the distribution of $\xi_1$ might be
specified in this way $P_x$ satisfies (EG). For instance, if $d=1$
and
$$
X_n=f(X_{n-1})+\xi_n,
$$
then for bounded $f$ and Laplacian random variable $\xi_1$ (EG)
holds. However, (EG) fails for many useful in applications ergodic
Markov chains. For $d=1$, a typical example gives Gaussian Markov
chain defined by a linear recurrent equation governed by i.i.d.
sequence of $(0,1)$-Gaussian random variables(here $|a|<1$)
$$
X_n=aX_{n-1}+\xi_n.
$$
In order to clarify this remark, notice that if (EG) holds true,
than for any bounded and measurable function $H$, satisfying
(H)-property, for some constants $K>0$, $\varrho\in (0,1)$, $n\ge
1,$
\begin{equation}\label{<const}
|E_xH(X_n)|\le K\varrho^n.
\end{equation}
However, the latter fails for $H(x)=\sign(x)$ satisfying
\eqref{nul}. In fact, if \eqref{<const} were correct, then $
\sum_{n=0}^\infty |E_xH(X_n)|\le \frac{K}{1-\varrho}. $ On the
other hand, it is readily to compute that $ \sum_{n=0}^\infty
|E_xH(X_n)| $ grows in $|x|$ on the set $\{|x|>1\}$ faster than
$O(\log(|x|)$.

\medskip
In this paper, we avoid a verification of (EG). Although our
approach is close to a conception of ``Multiplicative Ergodicity''
(see Balaji  and Myen \cite{BalMyen}) and ``Geometrical
Ergodicity'' (see Kontoyiannis and Meyn, \cite{KuMeyn} and Meyn
and Tweedie, \cite{metw}),  Chen and Guillin, \cite{ChenGui}) we
do not follow explicitly these methodologies.

Our main tools are the Poisson equation and the Puhalskii theorem
from \cite{P3}. The Poisson equation enables to reduce the MDP
verification for $(S^\alpha_n)_{n\ge 1}$ to
$(\frac{1}{n^\alpha}M_n)_{n\ge 1}$, where $M_n$ is a martingale
generated by Markov chain, while the Puhalskii theorem allows to
replace an asymptotic analysis for the Laplace transform of
$\frac{1}{n^\alpha}M_n$ by the asymptotic analysis for, so called,
{\it Stochastic Exponential}
\begin{equation}\label{StEx}
\mathscr{E}_n(\lambda)=\prod_{i=1}^nE\Big(\exp\Big[\Big\langle\lambda,
\frac{1}{n^\alpha} (M_i-M_{i-1})\Big\rangle\Big]\Big|X_{i-1}\Big),
\ \lambda\in\mathbb{R}^d
\end{equation}
being the product of the conditional Laplace transforms for
martingale increments.

An effectiveness of  the Poisson equation approach (method of
corrector) combined with the stochastic exponential is well known
from the proofs of functional central limit theorem (FCLT)  for
the family $(S^{0.5}_n)_{n\ge 1}$ (see, e.g. Papanicolaou, Stroock
and Varadhan \cite{PSV}, Ethier and Kurtz \cite{EK}, Bhattacharya
\cite{Br}, Pardoux and Veretennikov \cite{PV}; related topics can
be found in Metivier and Priouret (80's) for stochastic algorithms
analysis. The use of the same approach for a continuous time
setting can be found e.g. in \cite{LipSpo}, \cite{LSV}).

%=============================================================================

\section{Formulation of main result}
\label{sec+2} We consider Markov chain $(X_n)_{n\ge 0}$,
$X_n\in\mathbb{R}^d$ defined by a nonlinear recurrent equation
\begin{equation}\label{4.1a}
X_n=f(X_{n-1},\xi_n),
\end{equation}
where $f=f(z,v)$ is a vector function with entries
$f_1(z,v),\ldots,f_d(z,v)$, $u\in \mathbb{R}^d$, $v\in
\mathbb{R}^p$ and $(\xi_n)_{n\ge 1}$ is i.i.d. sequence of random
vectors of the size $p$.

We fix the following assumptions.

\medskip
\begin{assumption}\label{N.1}
Entries of $f$ are Lipschitz continuous functions in the following
sense: for any $v$
\begin{multline*}
|f_i(z_1\ldots,z_{j-1},z'_j,z_{j+1}\ldots,z_d,v_1,\ldots,v_p)
\\
- f_i(z_1\ldots,z_{j-1},z''_j,z_{j+1}\ldots,z_d,v_1,\ldots,v_p)|
\\
\le\varrho_{ij}|z'_j-z''_j|,
\\ \\
|f(z',v)-f(z'',v)|\le \varrho|z'-z''|,\\
\end{multline*}
where
$$
\max_{i,j}\varrho_{ij}=\varrho<1.
$$
\end{assumption}

\begin{assumption}\label{L.2}
For sufficiently small positive $\delta$, Cramer's condition
holds:
$$
Ee^{\delta|\xi_1|}<\infty.
$$
\end{assumption}

\medskip
\noindent
\begin{theorem}\label{theo-4.1}
Under Assumptions \ref{N.1} and \ref{L.2}, the Markov chain is
ergodic with the invariant measure $\mu$ such that
$\int_{\mathbb{R}^d}|z|\mu(dz)<\infty$.  For any Lipschitz
continuous function $H$, satisfying \eqref{nul}, the family
$(S^\alpha_n)_{n\ge 1}$ obeys the MDP in the metric space
$(\mathbb{R}^d,r)$ {\rm (}$r$ is the Euclidean metric{\rm )} with
the rate of speed $n^{-(2\alpha-1)}$ and the rate function given
in \eqref{I(x)2}.
\end{theorem}

\begin{remark}
Notice that:

- assumptions of Theorem \ref{theo-4.1} do not guarantee {\rm
(EG)},

- Lipschitz continuous $H$, obeying the linear growth condition,
are permis-

\hskip .1in sible for the MDP analysis,

- the $\xi_1$-distribution with a continuous component is not
required.
\end{remark}

\bigskip
Consider now a linear version of \eqref{4.1a}:
$$
X_n=AX_{n-1}+\xi_n,
$$
where $A$ is the $d\times d$-matrix with entries $A_{ij}$. Now,
Assumption \ref{N.1}  reads as: $\max_{ij}|A_{ij}|<1$. This
assumption is too restrictive. We replace it by more natural one

\begin{assumption}\label{L.1}
The eigenvalues of $A$ lie within the unit circle.
\end{assumption}

\begin{theorem}\label{theo-1}
Under Assumption \ref{L.1}, the Markov chain is ergodic with the
invariant measure $\mu$ such that
$\int_{\mathbb{R}^d}\|z\|^2\mu(dz)<\infty$. For any Lipschitz
continuous function $H$, satisfying \eqref{nul}, the family
$(S^\alpha_n)_{n\ge 1}$ obeys the MDP in the metric space
$(\mathbb{R}^d,r)$ with the rate of speed $n^{-(2\alpha-1)}$ and
the rate function given in \eqref{I(x)2}.
\end{theorem}

%=============================================================================
\section{Preliminaries}
\label{sec+3}
\subsection{(EG)-(H) conditions}
\mbox{} \label{sec+3.1} To clarify our approach to the MDP
analysis, let us first demonstrate  its applicability under
(EG)-(H) setting.

The (EG) condition provides the geometric ergodicity   of
$P^{(n)}_x$ to the invariant measure $\mu$ uniformly in $x$ in the
total variation norm: there exist constants $K>0$ and
$\varrho\in(0,1)$ such that for any $x\in\mathbb{R}^d$,
$$
\|P^{(n)}_x-\mu\|_{\sf tv}\le K\varrho^n, \ n\ge 1.
$$
The latter provides an existence of bounded function
\begin{equation}\label{ux}
U(x)=H(x)+\sum_{n\ge 1}P^{(n)}_xH
\end{equation}
solving the Poisson equation
\begin{equation}\label{Poisson1}
H(x)=H(x)+P_xU.
\end{equation}
In view of the Markov property, a sequence $(\zeta_i)_{i\ge 1}$ of
bounded random vectors with$ \zeta_i:=U(X_i)-P_{X_{i-1}}U$  forms
a martingale-differences relative to the filtration generated by
Markov chain. Hence, $M_n=\sum_{i=1}^n\zeta_i$ is the martingale
with bounded increments. With the help of Poisson's equation we
get the following decomposition
\begin{equation}\label{PoissDec}
\frac{1}{n^\alpha}\sum_{i=1}^nH(X_{i-1})=\underbrace{\frac{1}{n^\alpha}[U(x)-U(X_n)]}
_{\rm corrector}+\frac{1}{n^\alpha}M_n.
\end{equation}
The boundedness of $U$ provides a corrector negligibility in the
MDP scale, that is, the families $S^\alpha_n$ and
$\frac{1}{n^\alpha}M_n$ share the same MDP. In view of that,
suffice it to to establish the MDP for
$(\frac{1}{n^\alpha}M_n)_{n\ge 1}$.

\medskip
Assume for a moment that $\zeta_i$'s are i.i.d. sequence of random
vectors. Recall, $E\zeta_1=0$ and denote $B=E\zeta_1\zeta^*_1$.
Then, the Laplace transform for $\frac{1}{n^\alpha}M_n$ is:
\begin{equation}\label{Lap0}
\mathscr{E}_n(\lambda)=\Big(Ee^{\langle\lambda,
\frac{\zeta_1}{n^\alpha}\rangle}\Big)^n, \lambda\in\mathbb{R}^d.
\end{equation}
Under this setting, it is well known that $\frac{1}{n^\alpha}M_n$
obeys the MDP if $B$ is not singular matrix and
$$
\lim_{n\to\infty}n^{2\alpha-1}
\log\mathscr{E}_n(\lambda)=\frac{1}{2}\langle
\lambda,B\lambda\rangle, \ \lambda\in \mathbb{R}^d.
$$
We adapt this method of MDP verification to our setting. Instead
of $B$, we introduce matrices $ B(X_{i-1}), \ i\ge 1 $ with
\begin{equation}\label{B(x)}
B(x)=P_xUU^*-P_xU\big(P_xU)^*.
\end{equation}
The homogeneity of Markov chain and the definition of $\zeta_i$
provide a.s. that
$$
E(\zeta_i\zeta^*_i|X_{i-1})=B(X_{i-1}).
$$
Instead of the Laplace transform \eqref{Lap0}, we apply the
stochastic exponential \eqref{StEx}, expressed via $\zeta_i$'s,
$$
\mathscr{E}_n(\lambda)=\prod_{i=1}^nE\Big(e^{\langle\lambda,
\frac{\zeta_i}{n^\alpha}\rangle}\big|X_{i-1}\Big), \ \lambda\in
\mathbb{R},
$$
which is not the Laplace transform itself.

The Poisson equation \eqref{Poisson1} and its solution \eqref{ux}
permit to transform \eqref{B(x)} into
$$
B(x)=H(x)H^*(x)+\sum_{n\ge 1}\Big[H(x)\big(P^{(n)}_xH\big)^*+
\big(P^{(n)}_xH\big) H^* \Big],
$$
that is, $ \int_{\mathbb{R}^d}B(z)\mu(dz) $ coincides with $B$
from \eqref{B}.

\medskip
Now, we are in the position to formulate

\smallskip
\noindent {\bf Puhalskii Theorem.} [for more details, see
\cite{P3} and \cite{puh2}.] {\it Assume $B$ from \eqref{B} is
nonsingular matrix and for any $\varepsilon>0$,
$\lambda\in\mathbb{R}^d$
\begin{equation}\label{2.7Puh11}
\lim_{n\to\infty}\frac{1}{n^{2\alpha-1}}\log
P\Big(\Big|n^{2\alpha-1}
\log\mathscr{E}_n(\lambda)-\frac{1}{2}\langle\lambda,B\lambda\rangle\Big|>
\varepsilon\Big)=-\infty.
\end{equation}

Then, the family $\frac{1}{n^\alpha}M_n$, $n\ge 1$ possesses the
MDP in the metric space $(\mathbb{R}^d,r)$ {\rm (}$r$ is the
Euclidean metric{\rm )} with the rate of speed $n^{-(2\alpha-1)}$
and rate function $I(y)=\frac{1}{2}\|y\|^2_{B^{-1}}$. }

\medskip
\begin{remark}
The condition \eqref{2.7Puh11} is verifiable with the help of
\begin{equation}\label{I1d}
\begin{aligned}
& \lim_{n\to\infty}\frac{1}{n^{2\alpha-1}}\log P
\Big(\frac{1}{n}\Big|\sum_{i=1}^n\Big\langle\lambda,\big[B(X_{i-1})-B\big]\lambda
\Big\rangle \Big|>\varepsilon\Big)=-\infty
\\
& \lim_{n\to\infty}\frac{1}{n^{2\alpha-1}}\log P
\Big(\frac{1}{6n^{1+\alpha}}\sum_{i=1}^n
E\Big[|\zeta_i|^3e^{n^{-\alpha}|\zeta_i|}\big|X_{i-1}\Big]
>\varepsilon\Big)=-\infty.
\end{aligned}
\end{equation}
\end{remark}
The second condition in \eqref{I1d} is implied by the boundedness
of $|\zeta_i|$'s. The first part in \eqref{I1d} is known as
Dembo's conditions, \cite{Dem}, formulated as follows: for any
$\varepsilon>0$, $\lambda\in\mathbb{R}^d$
$$
\varlimsup_{n\to\infty}\frac{1}{n}\log P
\Big(\frac{1}{n}\Big|\sum_{i=1}^n\Big\langle\lambda,\big[B(X_{i-1})-B\big]\lambda
\Big\rangle \Big|>\varepsilon\Big)<0.
$$
In order to verify the first condition in \eqref{I1d}, we apply
again the Poisson equation technique. Set $
h(x)=\big\langle\lambda,\big[B(x)-B\big]\lambda\big\rangle $ and
notice that
$$
\int_{\mathbb{R}^d}h(z)\mu(dz)=0.
$$
Then, the function $ u(x)=h(x)+\sum_{n\ge 1}P^{(n)}_xh $ is well
defined and solves the Poisson equation $ u(x)=h(x)+P_xu. $
Similarly to \eqref{PoissDec}, we have
$$
\frac{1}{n}\sum_{i=1}^nh(X_{i-1})=\frac{u(x)-u(X_n)}{n}+
\frac{m_n}{n},
$$
where $m_n=\sum_{i=1}^nz_i$ is the martingale with bounded
martingale-differences $(z_i)_{i\ge 1}$. Since $u$ is bounded, the
first condition in \eqref{I1d} is reduced to
\begin{equation}\label{3.13a}
\lim_{n\to\infty}\frac{1}{n^{2\alpha-1}}\log
P\big(|m_n|>n\varepsilon\big)=-\infty
\end{equation}
while \eqref{3.13a} is provided by Theorem \ref{theo-A.1} in
Appendix which states that \eqref{3.13a} holds for any martingale
with bounded increments.

\subsubsection{Singular $B$}
\label{sec+3.1.1} The conditions from \eqref{I1d} remain to hold
whether $B$ is nonsingular or singular. For singular $B$ the
Puhalskii theorem is no longer valid. With singular $B$, we use
the Puhalskii theorem as helpful tool

It is well known that the family $\frac{M_n}{n^\alpha}$, $n\ge 1$
obeys the MDP with the rate of speed $n^{-(2\alpha-1)}$ and some
rate function,say $I(y)$ provided that
\begin{equation}\label{2.21sing1}
\begin{aligned}
& \varlimsup_{C\to \infty}\varlimsup_{n\to
\infty}\frac{1}{n^{2\alpha-1}}\log P
\Big(\Big\|\frac{M_n}{n^\alpha}\Big\|>C\Big)=-\infty
\\
& \varlimsup_{\varepsilon\to 0}\varlimsup_{n\to
\infty}\frac{1}{n^{2\alpha-1}}\log P
\Big(\Big\|\frac{M_n}{n^\alpha}-y\Big\|\le\varepsilon\Big)\le
-I(y)
\\
& \varliminf_{\varepsilon\to 0}\varliminf_{n\to
\infty}\frac{1}{n^{2\alpha-1}}\log P
\Big(\Big\|\frac{M_n}{n^\alpha}-y\Big\|\le \varepsilon\Big)\ge
-I(y).
\end{aligned}
\end{equation}
The first condition in \eqref{2.21sing1} provides the exponential
tightness in the metric $r$ while the next others the local MDP.

In order to verify of \eqref{2.21sing1}, we introduce
``regularized'' family $\frac{M^\beta _n}{n^\alpha},n\ge 1$ with
$$
M^\beta _n=M_n+\sqrt{\beta }\sum_{i=1}^n\vartheta_i,
$$
where $\beta $ is a positive parameter and $(\vartheta_i)_{i\ge
1}$ is a sequence of zero mean i.i.d. Gaussian random vectors with
$\cov(\vartheta_1,\vartheta_1)=:\mathbf{I}$ ($\mathbf{I}$ is the
unit matrix). The Markov chain and  $(\vartheta_i)_{i\ge 1}$ are
assumed to be independent objects.

It is clear that for this setting the matrix $B$ is transformed
into a positive definite matrix $B_\beta=B+\beta\mathbf{I}$. Now,
the Puhalskii theorem is applicable and guarantees the MDP with
the same rate of speed and the rate function
$$
I_\beta(y)=\frac{1}{2}\|y\|^2_{B^{-1}_ \beta }.
$$
We use now the well known fact (see, e.g. Puhalskii, \cite{P1})
that MDP provides the exponentially tightness and the the local
MDP:
\begin{equation}\label{2.21song1}
\begin{aligned}
& \varlimsup_{C\to \infty}\varlimsup_{n\to
\infty}\frac{1}{n^{2\alpha-1}}\log P \Big(\Big\|\frac{M^\beta
_n}{n^\alpha}\Big\|>C\Big)=-\infty
\\
& \varlimsup_{\varepsilon\to 0}\varlimsup_{n\to
\infty}\frac{1}{n^{2\alpha-1}}\log P \Big(\Big\|\frac{M^\beta
_n}{n^\alpha}-y\Big\|\le\varepsilon\Big)\le -I_\beta (y)
\\
& \varliminf_{\varepsilon\to 0}\varliminf_{n\to
\infty}\frac{1}{n^{2\alpha-1}}\log P \Big(\Big\|\frac{M^\beta
_n}{n^\alpha}-y\Big\|\le\varepsilon\Big)\ge -I_\beta (y).
\end{aligned}
\end{equation}
Notice now that \eqref{2.21sing1} is implied by \eqref{2.21song1}
if
\begin{equation}\label{2.23+1}
\lim_{\beta \to 0}I_\beta (y)=
\begin{cases}
\frac{1}{2}\|y\|^2_{B^\oplus}, & B^\oplus By=y
\\
\infty, & \text{otherwise}
\end{cases}
\end{equation}
and
\begin{equation}\label{2.24end1}
\lim_{\beta \to 0}\varlimsup_{n\to
\infty}\frac{1}{n^{2\alpha-1}}\log P \Big(\Big\|\frac{\sqrt{\beta
}}{n^\alpha}\sum_{i=1}^n\vartheta_i\Big\|> \eta\Big)=-\infty,
\quad \ \forall \ \eta>0.
\end{equation}

Let $T$ be an orthogonal matrix transforming $B$ to a diagonal
form: $ \diag(B)=T^*BT. $ Then, owing to
$$
2I_\beta (y)=y^*(\beta  I+B)^{-1}y=y^*T(\beta
I+\diag(B))^{-1}T^*y,
$$
for $y=B^\oplus By$ we have (recall that $B^\oplus
BB^\oplus=B^\oplus$, see \cite{Albert})
$$
\begin{aligned}
2I_\beta (y)&=y^*B^\oplus BT(\beta I+\diag(B))^{-1}T^*y
\\
&=y^*B^\oplus TT^* BT(\beta I+\diag(B))^{-1}T^*y
\\
&=y^*B^\oplus T\diag(B)(\beta I+\diag(B))^{-1}T^*y
\\
&\xrightarrow[\beta \to 0]{}y^*B^\oplus T\diag(B)\diag((B))^\oplus
T^*y
\\
&=y^*B^\oplus T\diag(B)T^*T(\diag(B))^\oplus T^*y
\\
&=y^*B^\oplus BB^\oplus u=u^*B^\oplus y=\|y\|^2_{B^\oplus}=2I(y).
\end{aligned}
$$

If $y\ne B^\oplus By$, $\lim_{\beta \to 0}2I_\beta (y)=\infty$.

Thus, \eqref{2.23+1} holds true.

\medskip
Since $(\vartheta_i)_{i\ge 1}$ is i.i.d. sequence of random
vectors and entries of $\vartheta_1$ are i.i.d. $(0,1)$-Gaussian
random variables, the verification of \eqref{2.24end1} is reduced
to
\begin{equation}\label{Ibeta1}
\lim_{\beta \to
0}\varlimsup_{n\to\infty}\frac{1}{n^{2\alpha-1}}\log P\Big(
\Big|\sum_{i=1}^n\xi_i\Big|>\frac{n^\alpha\eta}{\sqrt{\beta
}}\Big) =-\infty,
\end{equation}
where $(\xi_i)_{i\ge 1}$ is a sequence of i.i.d. $(0,1)$-Gaussian
random variables, and it suffices to consider the case ``+'' only.
By the Chernoff inequality with $\lambda>0$, we find that
$$
P\big(\sum_{i=1}^n\vartheta_i>\frac{n^\alpha\eta}{\sqrt{\beta
}}\big) \le \exp\Big(-\lambda\frac{n^\alpha\eta}{\sqrt{\beta
}}+n\frac{\lambda^2} {2}\Big)
$$
while the choice of $\lambda=\frac{n^\alpha\eta}{n\sqrt{\beta }}$
provides
$$
\frac{1}{n^{2\alpha-1}}\log P\Big(
\sum_{i=1}^n\eta_i>\frac{n^\alpha\eta}{\sqrt{\beta }}\Big)\le -
\frac{\eta^2}{2\beta }\xrightarrow[\beta \to 0]{}-\infty.
$$

%====================================================================

\subsection{Virtual scenario}
\label{sec+3.2} \mbox{}

- (EG)-(H) are not assumed

- the ergodicity of Markov chain is checked

- $H$ is chosen to hold \eqref{nul}.

\medskip
\noindent {\bf (1)} Let \eqref{ux} hold. Hence, the function $U$
solves the Poisson equation and the decomposition from
\eqref{PoissDec} is valid with $M_n=\sum_{i=1}^n\zeta_i$, where
$$
\zeta_i=u(X_i)-P_{X_{i-1}}u.
$$

Let
$$
\begin{aligned}
& E\zeta^*_i\zeta_i\le{\rm const.}
\\
&
E\Big[|\zeta_i|^3e^{n^{-\alpha}|\zeta_i|}\big|X_{i-1}\Big]\le{\rm
const.}
\end{aligned}
$$

\medskip
\noindent {\bf (2)} With $B(x)$ and $B$ are defined in
\eqref{B(x)} and \eqref{B} respectively, set
$$
h(x)=\Big\langle \lambda,\big[B(x)-B\big]\lambda\Big\rangle, \
\lambda\in\mathbb{R}^d.
$$

Let

(i) $u(x)=h(x)+\sum_{n\ge 1}P^{(n)}_xh$ is well defined

(ii) for $z_i=u(X_i)-P_{X_{i-1}}u$,
$$
\begin{aligned}
& Ez^2_i\le{\rm const.}
\\
& E\Big[|z_i|^3e^{n^{-\alpha}|z_i|}\big|X_{i-1}\Big]\le{\rm
const.}
\end{aligned}
$$

\medskip
\noindent {\bf (3)} For any $\varepsilon>0$, let
\begin{eqnarray*}
& \lim\limits_{n\to\infty}\frac{1}{n^{2\alpha-1}} \log
P\big(|U(X_n)|>n^\alpha\varepsilon\big)=-\infty
\\
& \lim\limits_{n\to\infty}\frac{1}{n^{2\alpha-1}} \log
P\big(|u(X_n)|>n^\alpha\varepsilon\big)=-\infty.
\end{eqnarray*}

\medskip
Notice that (EG)-(H) provide {\bf (1)}-{\bf (3)} and even if
(EG)-(H) fail, {\bf (1)}-{\bf (3)} may fulfill. Moreover, {\bf
(1)}-{\bf (3)} guarantee the validity for all steps of the proof
given in Section \ref{sec+3.1}.

Thus, an ergodic Markov chain, possessing  {\bf (1)}-{\bf (3)},
obeys the MDP.

\medskip
The proof of Theorems \ref{theo-4.1} and \ref{theo-1} follows this
scenario.

%========================================================================

\section{The proof of Theorem \ref{theo-4.1}}
\label{sec+4}

\subsection{Ergodic property}

\begin{lemma}\label{lem-4.1}
Under Assumption \ref{N.1}, $(X_n)_{n\ge 0}$ possesses the unique
probability invariant measure $\mu$ with
$\int_{\mathbb{R}^d}|z|\mu(dz)<\infty$.
\end{lemma}

\begin{proof}
Let $\nu$ be a probability measure on $\mathbb{R}^d$ with
$\int_{\mathbb{R}^d}|x|\nu(dx)<\infty$ and let a random vector
$X_0$, distributed in the accordance to $\nu$, is independent of
$(\xi_n)_{n\ge 1}$. We initialize the recursion, given in
\eqref{4.1a}, by $X_0$. Let now $X_n$ is generated by
\eqref{4.1a}. Then, $
\mu^n(dz)=\int_{\mathbb{R}^d}P^{(n)}_x(dz)\nu(dx) $ defines the
distribution of $X_n$.

We show that the family $(\mu^n)_{n\ge 1}$ is tight in the
Levy-Prohorov metric:
$$
\lim\limits_{k\to\infty}\varlimsup\limits_{n\to\infty}\mu^n(|z|>k)=0.
$$
By the Chebyshev inequality, $\mu^n(|z|>k)\le \frac{E|X_n|}{k}$.
The tightness follows from $ \sup_{n\ge 1}E|X_n|<\infty. $
Further, since By Assumption \ref{N.1},
$$
\begin{aligned}
|X_n|&=|f(0,\xi_n)+(f(X_{n-1},\xi_n)-f(0,\xi_n))|
\\
&\le |f(0,\xi_n)|+|f(X_{n-1},\xi_n)-f(0,\xi_n))|
\\
&\le |f(0,\xi_n)|+\varrho|X_{n-1}|
\\
&\le |f(0,0)|+\ell|\xi_n|+\varrho|X_{n-1}|,
\end{aligned}
$$
the sequence $(E|X_n|)_{n\ge 1}$ solves a recurrent inequality
$$
E|X_n|\le |f(0,0)|+\ell E|\xi_1|+\varrho E|X_{n-1}|
$$
subject to $E|X_0|=\int_{\mathbb{R}^d}|x|\nu(dx)(<\infty)$. Hence,
we find that for any $n\ge 1$,
\begin{equation*}
E|X_n|\le E|X_0|+\frac{|f(0,0)|+\ell E|\xi_1|} {1-\varrho}.
\end{equation*}

Thus, the family $\{\mu_n\}$ is tight, so that, by the Prohorov
theorem, $\{\mu^n\}$ contains further subsequence $\{\mu^{n'}\}$
converging, as $n'\nearrow\infty$,  in the Levy-Prohorov metric to
a limit $\mu$ being a probability measure on $\mathbb{R}^d$: for
any bounded and continuous function $g$ on $\mathbb{R}^d$
$$
\lim_{n'\to\infty}\int_{\mathbb{R}^d}g(z)\mu^{n'}(dz)=\int_{\mathbb{R}^d}g(z)\mu(dz).
$$
Thence, for $g(z)=L\wedge |z|$ and $L>0$, it holds
$$
\int_{\mathbb{R}^d}(L\wedge|z|)\mu(dz)=\lim_{n'\to\infty}E(L\wedge|X_{n'}|)
\le\varlimsup_{n\to\infty}E|X_n|<\infty
$$
and, by the monotone convergence theorem,
$$
\int_{\mathbb{R}^d}|z|\mu(dz)\le\varlimsup_{n\to\infty}E|X_n|<\infty.
$$
The $\mu$ is regarded now as a candidate to be the unique
invariant measure. So, we shall verify
\begin{equation*}
\int_{\mathbb{R}^d}g(x)\mu(dx)=\int_{\mathbb{R}^d}P_xg\mu(dx).
\end{equation*}
for any nonnegative, bounded and continuous function $g$. For
notational convenience, write $X^x_n$ and $X^\nu_n$, if $X_0=x$
and $X_0$ is distributed in the accordance with $\nu$. By
Assumption \ref{N.1},
$$
|X^x_n-X^\nu_n|\le\varrho|X^x_{n-1}-X^\nu_{n-1}|, \ n\ge 1,
$$
that is, $|X^x_n-X^\nu_n|$ converges to zero exponentially fast as
long as $n\to\infty$. For any $x\in\mathbb{R}^d$, the latter
provides $
\lim_{n'\to\infty}Eg(X^x_{n'})=\int_{\mathbb{R}^d}g(x)\mu(dx). $
Since the Markov chain is homogeneous, we also find that
$$
\lim_{n'\to\infty}Eg(X^x_{n'+1})=\int_{\mathbb{R}^d}g(z)\mu(dz).
$$
On the other hand, owing to $Eg(X^x_{n'+1})=EP_{X^x_{n'}}g$, the
above relation is nothing but
$$
\lim\limits_{n'\to\infty}EP_{X^x_{n'}}g=\int_{\mathbb{R}^d}g(z)\mu(dz).
$$
Finally, owing to $P_xg=Eg(f(x,\xi_1))$, the function $P_xg$ of
argument $x$ is bounded and continuous. Consequently, $
\lim\limits_{n'\to\infty}EP_{X^x_{n'}}g=\int_{\mathbb{R}^d}P_xg\mu(dx).
$

\smallskip
Assume $\mu'$ is another invariant probability measure, $\mu'\ne
\mu$. Then, taking $X^\mu_0$ and $X^{\mu'}_0$, distributed in the
accordance to $\mu$ and $\mu'$ respectively and independent of
$(\xi_n)_{n\ge 1}$, we get two stationary Markov chains
$(X^\mu_n)$ and $(X^{\mu'}_n)$ defined on the same probability
space as:
$$
\begin{aligned}
X^\mu_n=f(X^\mu_{n-1},\xi_n)
\\
X^{\mu'}_n=f(X^{\mu'}_{n-1},\xi_n).
\end{aligned}
$$
By Assumption \ref{N.1},
$|X^\mu_n-X^{\mu'}_n|\le\varrho|X^\mu_{n-1}-X^{\mu'}_{n-1}|$, i.e.
$\lim\limits_{n\to\infty}|X^\mu_n-X^{\mu'}_n|=0$. Recall that both
processes $X^\mu_n$ and $X^{\mu'}_n$ are stationary with the
marginal distributions $\mu$ and $\mu'$ respectively. Hence, for
any bounded and continuous function $g:\mathbb{R}^d\to\mathbb{R}$,
$$
\Big|\int_{\mathbb{R}^d}g(x)\mu(dx)-\int_{\mathbb{R}^d}g(x)\mu'(dx)\Big|
\le E|g(X_n^\mu)-g(X^{\mu'}_n)|\xrightarrow[n\to\infty]{}0,
$$
that is, $\mu=\mu'$.
\end{proof}

\subsection{The verification of (1)}
\label{sec+4.2} Let $K$ be the Lipschitz constant for $H$. Then
$|H(x)|\le |H(0)|+K|x|$ and
$\int_{\mathbb{R}^d}|H(z)|\mu(dz)<\infty.$ By \eqref{nul},
$EH(X^\mu_n)\equiv 0$. Then,
$$
\begin{aligned}
|EH(X^x_n)|&=|E(H(X^x_n)-H(X^\mu_n)|
\\
&\le K\varrho^nE|x-X^\mu_n|\le K(1+|x|)\varrho^n.
\end{aligned}
$$
Therefore, $\sum_{n\ge 1}|EH(X^x_n|\le\frac{K}{1-\varrho}(1+|x|)$.
Consequently, the function $U(x)$, given in \eqref{ux}, is well
defined and solves the Poisson equation.

Recall that $\zeta_i=U(X_i)-P_{X_{i-1}}U$.

\begin{lemma}\label{lem-2.1a}
The function $U(x)$ possesses the following properties:

{\rm 1)} $U(x)$ is Lipschitz continuous;

{\rm 2)} $P_x(UU^*)-P_xU(P_xU)^*$ is bounded and Lipschitz
continuous;

{\rm 3)} For sufficiently small $\delta>0$ and any $i\ge 1$
$$
E\Big(\big|U(X_i)-P_{\mbox{}_{X_{i-1}}}U\big|^3e^{\delta|U(X_i)-
P_{\mbox{}_{X_{i-1}}}U|}\big|X_{i-1}\Big)\le \text{\rm const.}
$$
\end{lemma}
\begin{proof}
1) Since  by Assumption \ref{N.1},
$$
|X^{x'}_n-X^{x''}_n|\le \varrho|X^{x'}_{n-1}-X^{x''}_{n-1}|, \
|X^{x'}_0-X^{x''}_0|\le |x'-x''|,
$$
we have
\begin{equation}\label{4.6pi}
\begin{aligned}
|U(x')-U(x'')|&\le |H(x')-H(x'')|+\sum_{n\ge
1}E|H(X^{x'}_n)-H(X^{x''}_n)|
\\
&\le \frac{K}{1-\varrho}|x'-x''|.
\end{aligned}
\end{equation}

\medskip
2) Recall (see \eqref{B(x)})
$$
P_x(UU^*)-P_xU(P_xU)^*=B(x)
$$
and denote $B_{pq}(x)$, $p,q=1,\ldots,d$ the entries of matrix
$B(x)$. Also, denote by $U_p(x)$, $p=1,\ldots,d$ the entries of
$U(x)$. Since $B(x)$ is nonnegative definite matrix, suffice it to
show only that $B_{pp}(x)$'s are bounded functions. Denote $F(z)$
the distribution function of $\xi_1$. Taking into the
consideration \eqref{4.6pi} and Assumption \ref{N.1}, we get
$$
\begin{aligned}
B_{pp}(x)&=E\Big(U_p\big(f(x,\xi_1)\big)-\int_{\mathbb{R}^d}U_p\big(f(x,z)\big)
dF(z)\Big)^2
\\
&\le
\frac{(K\ell)^2}{(1-\varrho)^2}E\Big|\int_{\mathbb{R}^d}|\xi_1-z|
dF(z)\Big|^2 \le
4\frac{(K\ell)^2}{(1-\varrho)^2}E|\xi_1|^2<\infty.
\end{aligned}
$$

The Lipschitz continuity of $B_{pq}(x)$ is proved similarly. Write
$$
B_{pq}(x')-B_{pq}(x'')=:ab-cd,
$$
where
$$
\begin{aligned}
&
a=E\Big(U_p\big(f(x',\xi_1)\big)-\int_{\mathbb{R}^d}U_q\big(f(x',z)\big)dF(z)\Big)
\\
&
b=E\Big(U_q\big(f(x',\xi_1)\big)-\int_{\mathbb{R}^d}U_q\big(f(x',z)dF(z)\Big)
\\
&
c=E\Big(U_p\big(f(x'',\xi_1)\big)-\int_{\mathbb{R}^d}U_q\big(f(x'',z)\big)dF(z)\Big)
\\
&
d=E\Big(U_q\big(f(x'',\xi_1)\big)-\int_{\mathbb{R}^d}U_q\big(f(x'',z)\big)dF(z)\Big).
\end{aligned}
$$
Now, applying $ab-cd=a(b-d)+d(a-c)$ and taking into account
\eqref{4.6pi} and Assumption \ref{N.1}, we find that $|a|,|d|\le
\frac{2K\ell}{1-\varrho}E|\xi_1| $ and so
$$
|B_{pq}(x')-B_{pq}(x'')|\le\frac{4K^2\ell\varrho}{(1-\varrho)^2}E|\xi_1||x'-x''|.
$$

\medskip
3) By \eqref{4.6pi} and Assumption \ref{N.1}
$$
|U(X_i)-P_{X_{i-1}}U|\le
\frac{K\ell}{1-\varrho}\big(E|\xi_1|+|\xi_i|\big).
$$
\end{proof}

\subsection{The verification of (2)}
\label{sec+4.3} The properties of $B(x)$ to be bounded and
Lipschitz continuous provide the same properties for
$$
h(x)=\big\langle \lambda,\big[B(x)-B\big]\lambda\big\rangle.
$$

Hence {\bf (2)} is provided by {\bf (1)}.

\subsection{The verification of (3)}
\label{sec+4.4} Since $U$ and $u$ are Lipschitz continuous, they
possess the linear growth condition, e.g., $ |U(x)|\le C(1+|x|), \
\exists C>0. $ So, {\bf (3)} is reduced to the verification of
\begin{equation}\label{2.10a}
\lim_{n\to\infty}\frac{1}{n^{2\alpha-1}}\log P\Big(\big|X_n\big|>
\varepsilon n^\alpha\Big)=-\infty, \ \varepsilon>0.
\end{equation}

Due to Assumption \ref{N.1}, we have
$$
\begin{aligned}
|X_n|&\le |f(X_{n-1},\xi_n)|\le |f(0,\xi_n)|+\varrho|X_{n-1}|
\\
&\le |f(0,0)| +\varrho|X_{n-1}|+\ell|\xi_n|.
\end{aligned}
$$
Iterating this inequality with $X_0=x$ we obtain
$$
\begin{aligned}
|X_n|&\le \varrho^n |x|+|f(0,0)|\sum_{j=1}^n\varrho^{n-j}+\ell
\sum_{j=1}^n\varrho^{n-j} |\xi_j|
\\
&\le
|x|+\frac{|f(0,0)|}{1-\varrho}+\ell\sum_{j=0}^{n-1}\varrho^j|\xi_{n-j}|.
\end{aligned}
$$
Hence, \eqref{2.10a} is reduced to
\begin{equation}\label{4.9m}
\lim_{n\to\infty}\frac{1}{n^{2\alpha-1}}\log
P\Big(\sum_{j=0}^{n-1}\varrho^j|\xi_{n-j}| \ge
n^\alpha\varepsilon\Big)=-\infty.
\end{equation}
We verify \eqref{4.9m} with the help of Chernoff's inequality:
with $\delta$, involving in Assumption \ref{L.2}, and
$\gamma=\frac{\delta}{1-\varrho}$
$$
\begin{aligned}
P\Big(\sum_{j=0}^{n-1}\varrho^j|\xi_{n-j}| \ge
n^\alpha\varepsilon\Big)\le e^{-n^\alpha
\gamma\varepsilon}Ee^{\sum_{j=0}^{n-1}\gamma\varrho^j|\xi_{n-j}|}.
\end{aligned}
$$
The i.i.d. property for $\xi_j$'s provides
$$
Ee^{\sum_{j=0}^{n-1}\gamma\varrho^j|\xi_{n-j}|}=
Ee^{\sum_{j=0}^{n-1}\gamma\varrho^j|\xi_1|}\le
Ee^{\sum_{j=0}^\infty\gamma\varrho^j|\xi_1|}=Ee^{\delta|\xi_1|}<\infty
$$
and we get
$$
\frac{1}{n^{2\alpha-1}}\log
P\Big(\sum_{j=0}^{n-1}\varrho^j|\xi_{n-j}| \ge
n^\alpha\varepsilon\Big) \le -n^{1-\alpha}\delta\varepsilon
+\frac{\log
Ee^{\delta|\xi_1|}}{n^{2\alpha-1}}\xrightarrow[n\to\infty]{}-\infty.
$$

\section{The proof of Theorem \ref{theo-1}}
\label{sec-5}

The proof of this theorem differs from the proof of Theorem
\ref{theo-4.1} only in some details concerning to (L.1). So, only
these parts of the proof are given below.

\subsection{Ergodic property and invariant measure}
\label{sec-5.1}

Introduce $(\widetilde{\xi}_n)_{n\ge 1}$ the independent copy of
$(\xi_n)_{n\ge 1}$. Owing to
$$
X_n=A^nx+\sum_{i=1}^nA^{n-i}\xi_i=A^nx+\sum_{i=0}^{n-1}A^i\xi_{n-i},
$$
we introduce
\begin{equation}\label{4.2x}
\widetilde{X}_n=A^nx+\sum_{i=0}^{n-1}A^i\widetilde{\xi}_i
\end{equation}
and notice that the i.i.d. property of $(\xi_i)_{i\ge 1}$ provides
$ (X_n)_{n\ge 0}\stackrel{\rm law}{=}(\widetilde{X}_n)_{n\ge 0}. $

By Assumption \ref{L.1}, $A^n\to 0$, $n\to\infty$, exponentially
fast. Particularly,
$$
\sum_{i=0}^\infty \trace\big(A^i\cov(\xi_1,\xi_1)
(A^i)^*\big)<\infty,
$$
so that $
\lim\limits_{n\to\infty}\widetilde{X}_n=\sum_{i=0}^\infty A^i
\widetilde{\xi}_i$ a.s. and in $\mathbb{L}^2$ norm.

Thus, the invariant measure $\mu$ is generated by the distribution
function of $\widetilde{X}_\infty$. In addition,
$E\|\widetilde{X}_\infty\|^2=\sum\limits_{i=0}^\infty
\trace\big(A^i\cov(\xi_1,\xi_1) (A^i)^*\big)$, so that
$$
\int_{\mathbb{R}^d}\|z\|^2\mu(dz)<\infty.
$$

\subsection{The verification of (1) and {\bf (2)}}
\label{sec-5.2}

Due to
$$
(X^{x'}_n-X^{x''}_n)=A(X^{x'}_{n-1}-X^{x''}_{n-1}),
$$
we have $(X^{x'}_n-X^{x''}_n)=A^n(x'-x'')$. Let us transform the
matrix $A$ into a Jordan form $ A = T J T^{-1} $ and notice that
$A^n=TJ^nT^{-1}$. It is well known that the maximal absolute value
of entries of $J^n$ is $n|\lambda|^n$, where $|\lambda|$ is the
maximal absolute value among eigenvalues of $A$. By Assumption
\ref{L.1}, $|\lambda|<1$. So, there exist $K>0$ and $\varrho<1$
such that $|\lambda|<\varrho$. Then, entries $A^n_{pq}$ of $A^n$
are evaluated as: $|A^n_{pq}|\le K\varrho^n$. Hence,
$|X^{x'}_n-X^{x''}_n|\le K\varrho^n|x'-x''|$, $n\ge 1$, and the
verification of {\bf (1), (2)} is in the framework of Section
\ref{sec+3}.

\subsection{The verification of (3)}
\label{sec-5.3} As in Section \ref{sec+3}, the verification of
this property is reduced to
\begin{equation}\label{2.10aa}
\lim_{n\to\infty}\frac{1}{n^{2\alpha-1}}\log P\Big(\big|X_n\big|>
\varepsilon n^\alpha\Big)=-\infty, \ \varepsilon>0.
\end{equation}
In \eqref{2.10aa}, we may replace $X_n$ by its copy
$\widetilde{X}_n$ defined in \eqref{4.2x}. Notice also that
$$
|\widetilde{X}_n|\le |A^nx|+\sum_{i=0}^\infty
\max_{pq}|A^i_{pq}||\widetilde{\xi}|.
$$
As was mentioned above, $|A^i_{pq}|\le K\varrho^j$ for some $K>0$
and $\varrho\in(0,1)$. Hence, suffice it to verify
$$
\lim_{n\to\infty}\frac{1}{n^{2\alpha-1}}\log
P\Big(\sum_{i=0}^\infty \varrho^i |\xi_i|>\varepsilon
n^\alpha\Big)=-\infty, \ \varepsilon>0
$$
what be going on similarly to corresponding part of the proof in
Section \ref{sec+3}.

\section{Exotic example}
\label{sec-2.3} Let $(X_n)_{n\ge 0}$, $X_n\in\mathbb{R}$ and
$X_0=x$, be Markov chain defined by the recurrent equation
\begin{equation}\label{exmp1}
X_n=X_{n-1}-m\frac{X_{n-1}}{|X_{n-1}|}+\xi_n,
\end{equation}
where $m$ is a positive parameter, $(\xi_n)$ is i.i.d. sequence of
zero mean random variables with
$$
Ee^{\delta|\xi_1|}<\infty,  \ \text{for some $\delta>0$},
$$
and let $\frac{0}{0}=0$.

Although the virtual scenario is not completely verifiable here we
show that for
$$
H(x)=\frac{x}{|x|}
$$
the family $(S^\alpha_n)_{n\ge 1}$ possesses the MDP provided that
\begin{equation}\label{mdel}
m>\frac{1}{\delta}\log Ee^{\delta|\xi_1|}.
\end{equation}

Indeed, by \eqref{exmp1} we have
$$
\frac{1}{n^\alpha}\sum_{k=1}^n\frac{X_{k-1}}{|X_{k-1}|}=
\frac{1}{m}\frac{(X_n-x)}{n^\alpha}+\frac{1}{n^\alpha}
\sum_{k=1}^n\frac{\xi_k}{m}.
$$
The family
$\big(\frac{1}{n^\alpha}\sum\limits_{k=1}^n\frac{\xi_k}{m}\big)_{n\ge
1}$ possesses the MDP with the rate of speed $n^{-(2\alpha-1)}$
and the rate function $I(y)=\frac{m^2}{2E\xi^2_1}y^2$. Then, the
family $(S^\alpha_n)_{n\ge 1}$ obeys the same MDP provided that
$\big(\frac{X_n-x}{n^\alpha}\big)_{n\ge 1}$ is exponentially
negligible family with the rate $n^{-(2\alpha-1)}$. This
verification is reduced to
\begin{equation}\label{negl}
\lim_{n\to\infty}\frac{1}{n^{2\alpha-1}}\log
P\big(|X_n|>n^\alpha\varepsilon\big)= -\infty, \ \varepsilon>0.
\end{equation}
By the Chernoff inequality $
P\big(|X_n|>n^\alpha\varepsilon\big)\le e^{-\delta
n^\alpha\varepsilon} Ee^{\delta|X_n|}, $ that is \eqref{negl}
holds if $\sup\limits_{n\ge 1}Ee^{\delta|X_n|}<\infty$ for some
$\delta>0$. We show that the latter holds true for $\delta$
involved in \eqref{mdel}. A helpful tool for this verification is
the inequality $\big|z-m\frac{z}{|z|}\big|\le \big||z|-m\big|$.
Write
$$
\begin{aligned}
Ee^{\delta|X_n|}&=Ee^{\delta|X_n|}I(|X_{n-1}|\le
m)+Ee^{\delta|X_n|}I(|X_{n-1}|>m)
\\
&\le e^{\delta m}Ee^{\delta|\xi_1|}+e^{-\delta
m}Ee^{\delta|\xi_1|}Ee^{\delta|X_{n-1}|}.
\end{aligned}
$$
Set $\ell=e^{\delta m}Ee^{\delta|\xi_1|}$ and $\varrho=e^{-\delta
m}Ee^{\delta|\xi_1|}$. By \eqref{mdel}, $\varrho<1$. Hence,
$V(x)=e^{\delta|x|}$ is the Lyapunov function: $ P_xV\le\varrho
V(x)+\ell. $ Consequently,
$$
EV(X_n)\le \varrho EV(X_n)+\ell, \  n\ge 1
$$
and so, $\sup_{n\ge 1}EV(X_n)\le V(x)+\frac{\ell}{1-\varrho}$.

\section{Statistical example}
\label{sec-7}

An asymptotic analysis, given in this section, demonstrate the
thesis ``MDP instead of CLT''.

Let
$$
X_n=\theta f(X_{n-1})+\xi_n,
$$
where $\theta$ is a number and $(\xi_n)_{n\ge 1}$ is i.i.d.
sequence of of $(0,1)$-Gaussian random variables. We assume that
$|\theta|<1$ and $f$ is bounded continuously differentiable
function with $|f'(x)|\le 1$.
By Theorem \ref{theo-4.1}, $(X_n)$
is an ergodic Markov chain and its invariant measure $\mu_\theta$
depends on parameter $\theta$. Since $\xi_1$ is Gaussian random variables,
$\mu_\theta$, being a convolution of some measure with Gaussian one, possesses
a density relative to $dz$. Then, assuming $f^2(x)>0$ relative to Lebesgue measure,
we have
$ B_\theta=\int_\mathbb{R}f^2(z)\mu(dz)>0.$ Under the above assumptions,
\begin{equation*}
\theta_n=\frac{\sum_{i=1}^nf(X_{i-1})X_i}{\sum_{i=1}^nf^2(X_{i-1})}
\end{equation*}
is a strongly consistent estimate of $\theta$ by sampling $\{X_1,\ldots,X_n\}$,
that is, $\lim_{n\to\infty}\theta_n=\theta$ a.s. Moreover, it is known its asymptotic
in the CLT scale:
$$
\sqrt{n}(\theta-\theta_n)\xrightarrow[n\to\infty]{\rm
law}\Big(0,\frac{1}{B_\theta} \Big)\text{-Gaussian r. v.}
$$
Here, we give an asymptotic of $\theta_n$ in the MDP scale:
for any $\alpha\in \big(\frac{1}{2},1\big)$,
$$
n^{1-\alpha}(\theta-\theta_n)\xrightarrow[n\to\infty]{\rm MDP}
\Big(\frac{1}{n^{2\alpha-1}},\frac{y^2}{2B_\theta}\Big).
$$
\begin{theorem}\label{7.1c}
The family $ n^{1-\alpha}(\theta-\theta_n) $ obeys the MDP with
the rate of speed $\frac{1}{n^{2\alpha-1}}$ and the rate function
$I(y)=\frac{y^2}{2B_\theta}$.
\end{theorem}
\begin{proof}
The use of
$$
n^{1-\alpha}(\theta-\theta_n)
=\frac{\frac{1}{n^\alpha}\sum_{i=1}^nf(X_{i-1})\xi_i}{\frac{1}{n}\sum_{i=1}^n
f^2(X_{i-1})}
$$
and the law of large numbers,
$P\text{-}\lim_{n\to\infty}\frac{1}{n}\sum_{i=1}^n
f^2(X_{i-1})=B_\theta$, give a hint that that the theorem
statement is valid provided that

(i) for $M_n=\sum_{i=1} f(X_{i-1})\xi_i$,
the family $\big(\frac{1}{n^\alpha}M_n\big)_{n\to\infty}$ obeys the MDP

\hskip .17in with the rate of speed $\frac{1}{n^{2\alpha-1}}$ and
the rate function $I(y)=\frac{y^2}{2B^{-1}_\theta}$;

(ii) for any $\varepsilon>0$,
$$
\lim_{n\to\infty}\frac{1}{n^{2\alpha-1}}\log
P\Big(\Big|\frac{1}{n}\sum_{i=1}^n\Big[
f^2(X_{i-1})-B_\theta\Big]\Big|\ge\varepsilon\Big)=-\infty.
$$
Following to \eqref{StEx} and taking into account the setting, we
notice that
$$
\mathscr{E}_n(\lambda)=\exp\Big(\sum_{i=1}^n
\frac{\lambda^2}{2n^{2\alpha}} f^2(X_{i-1})\Big).
$$
is the stochastic exponential related to
$\big(\frac{1}{n^\alpha}M_n\big)_{n\to\infty}$.
Consequently, \eqref{2.7Puh11} is reduced to (ii), that is, only
(ii) is left to be verified.

The verification of (ii) is in the framework of Theorem \eqref{theo-4.1}.
The function $H(x)=f^2(x)-B_\theta$ satisfies the assumptions of Theorem
\ref{theo-4.1}. Hence, the family
$\big(\frac{1}{n^\alpha}\sum _{i=k}^nH(X_{k_1})\big)_{n\to\infty}$
obeys the MDP with the rate of speed $\frac{1}{n^{2\alpha-1}}$ and
the rate function
$$
J(y)=
  \begin{cases}
\frac{y^2}{2}\widehat{B}^\oplus_\theta & \widehat{B}_\theta>0, \\
\infty &, \widehat{B}_\theta=0, \ y\ne 0,
  \end{cases}
$$
where, in accordance with \eqref{B},
$$
\widehat{B}_\theta=\int_\mathbb{R}H^2(x)\mu_\theta(dx) +2\sum_{n\ge
1}\int_\mathbb{R} H(x)P^{(n)}_xH\mu_\theta(dx).
$$
In particular,
$$
\varlimsup_{n\to\infty}\frac{1}{n^{2\alpha-1}}\log
P\Big(\Big|\frac{1}{n^\alpha}\sum_{k=1}^nH(X_{k-1})\Big| \ge
C\varepsilon\Big)\le
  \begin{cases}
    -\frac{1}{2\widehat{B}_\theta}C^2\varepsilon^2, & \widehat{B}^\theta>0 \\
    -\infty, & \text{otherwise}.
  \end{cases}
$$
Hence, for any $C>0$, we find that
\begin{equation*}
\begin{aligned}
& \varlimsup_{n\to\infty}\frac{1}{n^{2\alpha-1}}\log
P\Big(\Big|\frac{1}{n}\sum_{k=1}^nH(X_{k-1})\Big|
\ge\varepsilon\Big)
\\
&=\varlimsup_{n\to\infty}\frac{1}{n^{2\alpha-1}}\log
P\Big(\Big|\frac{1}{n^\alpha}\sum_{k=1}^nH(X_{k-1})\Big| \ge
n^{1-\alpha}\varepsilon\Big)
\\
&\le \varlimsup_{n\to\infty}\frac{1}{n^{2\alpha-1}}\log
P\Big(\Big|\frac{1}{n^\alpha}\sum_{k=1}^nH(X_{k-1})\Big| \ge
C\varepsilon\Big)
\\
&\le
  \begin{cases}
-\frac{C^2\varepsilon^2}{2\widehat{B}_\theta}     & \widehat{B}_\theta>0, \\
    -\infty & \text{otherwise}
  \end{cases}
  \xrightarrow[C\to\infty]{}-\infty.
\end{aligned}
\end{equation*}
\end{proof}

\appendix
\section{Exponentially integrable martingale-differences}

Let $\zeta_n=(\zeta_n)_{n\ge 1}$ be a martingale-difference  with
respect to some filtration $\mathscr{F}=(\mathscr{F}_n)_{n\ge 0}$
and $M_n=\sum\limits_{i=1}^n\zeta_i$ be the corresponding
martingale.

\begin{theorem}\label{theo-A.1}
Assume that for sufficiently small positive $\delta$ and any $i\ge
1$
\begin{equation}\label{A.1}
E\big(e^{\delta|\zeta_i|}|\mathscr{F}_{i-1}\big)\le {\rm const.}
\end{equation}
Then for any $\alpha\in (0.5,1)$
$$
\lim_{n\to\infty}\frac{1}{n^{2\alpha-1}}\log
P\big(|M_n|>n\varepsilon\big)=-\infty.
$$
\end{theorem}
\begin{proof}
Suffice it to prove $
\lim\limits_{n\to\infty}\frac{1}{n^{2\alpha-1}}\log P\big(\pm
M'_n>n \varepsilon\big)=-\infty$. We verify here only ``+'' only
(the proof of ``-'' is similar).

\medskip
For fixed positive $\lambda$ and sufficiently large $n$, let us
introduce the stochastic exponential
$$
\mathscr{E}_n(\lambda)=\prod_{i=1}^nE\big(e^{\lambda
\frac{\zeta_i}{n}}\big|\mathscr{F}_{i-1}\big).
$$
A direct verification shows that
$$
E\exp\Big(\frac{\lambda M_n}{n}-\log
\mathscr{E}_n(\lambda)\Big)=1.
$$
We apply this equality for further ones
\begin{equation}\label{StTv}
\begin{aligned}
1&\ge EI\Big(M_n>n\varepsilon\Big) \exp\Big(\frac{\lambda
M_n}{n}-\log  \mathscr{E}_n(\lambda)\Big)
\\
&\ge EI\Big(M_n>n\varepsilon\Big) \exp\Big(\lambda\varepsilon-\log
\mathscr{E}_n(\lambda)\Big).
\end{aligned}
\end{equation}
Due to $E\big(\lambda\frac{\zeta_i}{n}|\mathscr{F}_{i-1}\big)=0$
and \eqref{A.1}, we find that
$$
\begin{aligned}
\log  \mathscr{E}_n(\lambda)&=\sum_{i=1}^n\log\Big(1+E\big[
e^{\lambda\frac{\zeta_i}{n}}-1-\lambda\frac{\zeta_i}{n}|\mathscr{F}_{i-1}\big]\Big)
\\
&
\le\sum_{i=1}^n\Big\{\frac{\lambda^2}{2n^2}E\big((\zeta_i)^2|X_{i-1}\big)+
\frac{\lambda^3}{6n^3}E\big(|\zeta_i|^3e^{\lambda\frac{|\zeta_i|}{n}}|\mathscr{F}
_{i-1}\big) \Big\}
\\
& \le K\Big[\frac{\lambda^2}{2n}+\frac{\lambda^3}{6n^2}\Big],
\end{aligned}
$$
where $K$ is some constant. This inequality, being incorporated
into \eqref{StTv}, provides
$$
1\ge EI\big(M_n>n\varepsilon\big)
\exp\Big(\lambda\varepsilon-K\Big[\frac{\lambda^2}{2n}+\frac{\lambda^3}{6n^2}
\Big]\Big).
$$
If $\varepsilon<3$, taking $\lambda=\varepsilon nK^{-1}$, we find
that
$$
\frac{1}{n^{2\alpha-1}}\log P\big(M_n>n\varepsilon\big)\le
-\frac{\varepsilon^2n^{2(1-\alpha)}}{K}\Big(\frac{1}{2}-\frac{\varepsilon}{6}\Big)
\xrightarrow[n\to\infty]{}-\infty
$$

Thus, the desired statement holds true.
\end{proof}

\end{document}